\documentclass{notemata}
\usepackage{url}
\usepackage{amsmath}
\usepackage{amssymb}
\usepackage{nmmacro}
\chardef\bslchar=`\\ 
\makeatletter
\newcommand{\addbslash}{\expandafter\@addbslash\string}
\def\@addbslash#1{\bslchar\@nobslash#1}
\newcommand{\nobslash}{\expandafter\@nobslash\string}
\def\@nobslash#1{\ifnum`#1=\bslchar\else#1\fi}
\newcommand{\ntt}{\normalfont\ttfamily}
\def\@boxorbreak{\leavevmode
  \ifmmode\hbox\else\ifdim\lastskip=\z@\penalty9999 \fi\fi}
\DeclareRobustCommand{\cs}[1]{\@boxorbreak{\ntt\addbslash#1\@empty}}
\makeatother
\begin{document}
\begin{article}
\begin{opening}
  \title{Inaccuracy measures for concomitants of generalized order statistics in Morgenstern family}
\author{Safieh Daneshi
\thanks{This work is partially supported by the Clinical Research Development Center (CRDC) of the Bushehr University of Medical Sciences}}
  \institute{Department of Statistics, Faculty of Mathematical Sciences, Shahrood University of Technology, Shahrood, Iran.\\
    \email{s.daneshi@bpums.ac.ir}}
\author{Ahmad Nezakati}
 \institute{Department of Statistics, Faculty of Mathematical Sciences, Shahrood University of Technology, Shahrood, Iran.\\
    \email{nezakati@shahroodut.ac.ir }}
  \author{Saeid Tahmasebi}
  \institute{Department of Statistics, Persian Gulf University, Bushehr, Iran.\\
    \email{tahmasebi@pgu.ac.ir}}
  \author{Maria Longobardi  \thanks{This work is partially supported by the GNAMPA research group of INDAM (Istituto Nazionale di Alta Matematica) and MIUR-PRIN 2017}}
  \institute{Dipartimento di Matematica e Applicazioni Universit$\acute{a}$ di Napoli Federico II Via Cintia, I-80126 Napoli, Italy.\\
    \email{maria.longobardi@unina.it}}
  \runningauthor{S. Daneshi, A. Nezakati, S.Tahmasebi, M.Longobardi} \runningtitle{Inaccuracy measures for concomitants of GOS in Morgenstern family}

\begin{abstract}
 In this paper, we obtain a measure of inaccuracy between rth concomitant of generalized order statistic and the parent random variable in Morgenstern family. Applications of this result are given for concomitants of order statistics and record values. We also study some results of cumulative past inaccuracy (CPI) between the distribution function of rth concomitant of order statistic (record value) and the distribution function of parent random variable. Finally, we discuss on a problem of estimating the CPI by means of the empirical CPI in concomitants of generalized order statistics.
\end{abstract}

\end{opening}

\section{Introduction}

Let $(X_{i},Y_{i})$,  $i=1,2,\cdots,n$  be independent and identically distributed random variables from a continuous bivariate distribution $F_{X,Y}(x,y)$. If $X_{(r:n)}$ denotes the $r$th order statistic, then the $Y$'s  associated with $X_{(r:n)}$ denoted by $Y_{[r:n]}$ is called the concomitant of $r$th order statistic. The concomitants are of interest in selection and prediction problems.
The concept of generalized order statistics (GOS) was introduced by [7] as a unified approach to a variety of models of ordered random variables such as ordinary order statistics, sequential order statistics, progressive type-II censoring, record values and Pfeifers records. The random variables $X(1, n,m, k), X(2, n,m, k),\cdots,X(n, n,m, k)$ are called
generalized order statistics based on the absolutely continuous distribution function(cdf) $F$ with density function $f$, if their joint density function is given by

\begin{eqnarray*}
‎f^{X(1,n,m,k),...,X(n,n,m,k)}(x_1,...,x_n)&=&k\left(\prod_{j=1}^{n-1}\gamma_j\right)\left(\prod_{i=1}^{n-1}(1-F(x_i))^m f(x_i)\right)\nonumber \\
&\times&(1-F(x_n))^{k-1}f(x_n),\\‎
‎&&F^{-1}(0)\leq x_1\leq x_2\leq‎ ... ‎\leq x_n \leq F^{-1}(1)‎,
\end{eqnarray*}

with parameters $n\in \mathbb{N}, k>0,m\in\mathbb{R}$, such that $\gamma_r=k+(n-r)(m+1)>0$, for all $1\leq r\leq n$. Similarly, concomitants can also be defined in the case of GOS.

[9] defined a class of bivariate distributions with the probability density function (pdf) given by
\begin{equation}{\label{Mf}}
f_{X,Y}(x,y)=f_{X}(x)f_{Y}(y)\left[1+\alpha (2F_{X}(x)-1)(2F_{Y}(y)-1)\right], \;  \;  \; \;  |\alpha|\leq 1,
 \end{equation}
where $\alpha$ is the association parameter (see [5], and references therein for more details). For the Morgenstern family with pdf given by \eqref{Mf}, the density function and distribution function of the concomitant of rth GOS (denoted by $Y_{[r,n,m,k]}, 1\leq r\leq n$), are given by [2] as follows:
\begin{equation}{\label{g1}}
g_{[r,n,m,k]}(y)=f_Y(y)\left[1+\alpha C^*(r,n,m,k)(1-2F_Y(y))\right],
\end{equation}

\begin{equation}{\label{G1}}
G_{[r,n,m,k]}(y)=F_Y(y)\left[1+\alpha C^*(r,n,m,k)(1-F_Y(y))\right],
\end{equation}
where $C^*(r,n,m,k)=\frac{2\prod_{j=1}^r \gamma_j}{\prod_{i=1}^r (\gamma_i+1)}-1$. In the special case of GOS in Morgenstern family, if $Y_{[r:n]}$ denotes the concomitant of rth order statistic $X_{(r:n)}$, then the pdf and cdf of $Y_{[r:n]}$ in Morgenstern family are given by
\begin{eqnarray*}
  f_{Y_{[r:n]}}(y)=f_{Y}(y)\left[1+\alpha\left(\frac{n-2r+1}{n+1}\right)(1-2F_{Y}(y))\right],
\end{eqnarray*}
and
\begin{eqnarray*}
  F_{Y_{[r:n]}}(y)=F_{Y}(y)\left[1+\alpha\left(\frac{n-2r+1}{n+1}\right)(1-F_{Y}(y))\right],
\end{eqnarray*}
 respectively. We refer the reader to [1] for more details.

Let $(X_{1},Y_{1}),(X_{2},Y_{2}),\cdots$ be a sequence of bivariate random
variables from a continuous distribution. If $\{R_n,n\geq1\}$ is the sequence of upper record values in the sequence of $X$'s, then the $Y$ which corresponds with
the nth-record will be called the concomitant of the nth-record, denoted by $
R_{[n]}$. The concomitants of record values arise in a wide variety of practical experiments such as industrial stress testing, life time experiments, meteorological analysis, sporting matches and some other experimental fields. For other important applications of record values and
their concomitants see [1]. The pdf and cdf for $R_{[n]}$ has been obtained as follows:
\begin{equation}
f_{R_{[n]}}(y)=f_Y(y)[1+%
\alpha_{n}(1-2F_Y(y))],\;\; n\geq 1,
\end{equation}
and
\begin{equation}
F_{R_{[n]}}(y)=F_Y(y)[1+%
\alpha_{n}(1-F_Y(y))],
\end{equation}
where $\alpha_{n}=\alpha(2^{1-n}-1)$ .

Let $X$ and $Y$ be two non-negative random variables with distribution functions $F(x)$ and $G(x)$, respectively. If $f (x)$ is the actual probability density function (pdf) corresponding to the observations and $g(x)$ is the density assigned by the experimenter, then the inaccuracy measure of $X$ and $Y$ is defined by [8] as follows:
 \begin{equation*}
I(f,g)=-\int_{0}^{+\infty}f(x)\log g(x)dx.
\end{equation*}

Analogous to this measure of inaccuracy, [12] proposed a cumulative past inaccuracy (CPI) measure as
\begin{equation*}
I(F,G)=-\int_{0}^{+\infty}F(x)\log G(x)dx.
\end{equation*}

 Several authors have worked on measures of inaccuracy  for ordered random variables. [11] proposed the measure of inaccuracy between the ith order statistic and the parent random variable. [12] developed measures of dynamic cumulative residual and past inaccuracy. They studied characterization results of these dynamic measures under proportional hazard model and proportional reversed hazard model. Recently [13] have introduced the measure of residual inaccuracy of order statistics and proved a characterization result for it. Motivated by some of the articles mentioned above, in this paper we aim to present some results on inaccuracy for concomitants of GOS in Morgenstern family. The paper is organized as follows: In Section 2, we obtain a measure of inaccuracy between $g_{[r,n,m,k]}(y)$ and $f_Y(y)$ in Morgenstern family. Another applications of this result are given for concomitants of order statistics and record values. We also study some results of CPI  between $G_{[r,n,m,k]}(y)$ and $F_Y(y)$. Applications of CPI are given for concomitants of order statistics and record values. Finally, in Section 3, we discuss on a problem of estimating the CPI by means of the empirical CPI for concomitants of GOS.

\section{Inaccuracy measures for concomitants of GOS}
If $Y_{[r,n,m,k]}$ is the concomitant of rth GOS from \eqref{Mf},
then the inaccuracy measure between $g_{[r,n,m,k]}(y)$ and $f_{Y}(y)$ for $1\leq r \leq n$,  $\alpha\neq 0$ is given by
\begin{eqnarray}\label{Ig1}
  I(g_{[r,n,m,k]},f_{Y})&=&-\int_0^{\infty}g_{[r,n,m,k]}(y)\log{f_{Y}(y)}dy \nonumber\\
  &=&\left[1+\alpha C^*(r,n,m,k)\right]H(Y)\nonumber \\
  &+&2\alpha C^*(r,n,m,k)\int_0^{\infty}f_{Y}(y)F_{Y}(y)\log{f_{Y}(y)}dy\nonumber\\
  &=&\left[1+\alpha C^*(r,n,m,k)\right]H(Y)\nonumber \\
  &+&2\alpha C^*(r,n,m,k)\int_0^{1}u\log f_Y(F^{-1}_Y(u))du \nonumber\\
 &=&\left[1+\alpha C^*(r,n,m,k)\right]H(Y)+2\alpha C^*(r,n,m,k)\phi_f(u),
\end{eqnarray}
where $\phi_f(u)=\int_0^{1}u\log f_Y(F_Y^{-1}(u))du$ and  $$H(Y)=-\int_0^{\infty}f_{Y}(y)\log{f_{Y}(y)}dy$$
 is the Shannon entropy  of the random variable $Y$.

As an application of the representation \eqref{Ig1}, we consider the following special cases.\\
\textbf{Case 1:} According to \eqref{Ig1}, if we put $m = 0$ and $k = 1$, then an inaccuracy measure  between $f_{Y_{[r:n]}}$ (density function of rth concomitant of order statistic) and $f_{Y}$ in Morgenstern family is obtained as follows:
\begin{eqnarray}\label{Ig2}
I(f_{Y_{[r:n]}},f_{Y})&=&-\int_0^{\infty}f_{Y_{[r:n]}}(y)\log{f_{Y}(y)}dy \nonumber\\
    &=&\left[1+\alpha \left(\frac{n-2r+1}{n+1}\right)\right] H(Y) +2\alpha \left(\frac{n-2r+1}{n+1}\right)\phi_f(u)\nonumber\\
  &=&H(Y)+\frac{3(n-2r+1)}{2(n+1)}\left[I(f_{Y_{[1:2]}},f_{Y})-I(f_{Y_{[2:2]}},f_{Y})\right].
  \end{eqnarray}
In the following, we present some examples and properties of $I(f_{Y_{[r:n]}},f_{Y})$.
\begin{Example}
Let $(X_{i},Y_{i})$, $i=1,2,...,n$ be a random sample from Gumbel
bivariate exponential distribution (GBED) with cdf
\begin{equation}
F(x,y)=\left(1-\exp\left(\frac{-x}{\theta_{1}}\right)\right)\left(1-\exp\left(\frac{-y}{\theta_{2}}\right)\right)\left[1+\alpha\exp\left(\frac{-x}{\theta_{1}}-\frac{y}{\theta_{2}}\right)\right].
\end{equation}
From \eqref{Ig2}, we find
\begin{eqnarray}
\label{eq2}
  I(f_{Y_{[r:n]}},f_{Y})=[1+\log \theta_{2}]-\frac{\alpha}{2}\left(\frac{n-2r+1}{n+1}\right).
\end{eqnarray}
By using \eqref{eq2}, we get
\begin{eqnarray*}
 A_{\alpha}(n)= I(f_{Y_{[n:n]}},f_{Y})- I(f_{Y_{[1:n]}},f_{Y})=\alpha\left(\frac{n-1}{n+1}\right),
\end{eqnarray*}
which is positive, negative or zero whenever $0 <\alpha\leq1, n > 1$; $-1\leq\alpha < 0, n > 1$ or
$n = 1 \; or\; \alpha = 0$, respectively. Also, the difference between $I(f_{Y_{[r:n]}},f_{Y})$ and $H(Y)$ is
\begin{eqnarray*}
 B_{\alpha,n}(r)=I(f_{Y_{[r:n]}},f_{Y})- H(Y)=-\frac{\alpha}{2}\left(\frac{n-2r+1}{n+1}\right).
\end{eqnarray*}
$B_{\alpha,n}(r)$ is positive for $-1\leq\alpha<0$ , $1\leq r<\frac{n+1}{2}$ (or $0<\alpha\leq1$, $\frac{n+1}{2}< r\leq n)$. Also, it is negative for
  $-1\leq\alpha<0$ , $\frac{n+1}{2}< r\leq n$( or $0<\alpha\leq1$ ,$1\leq r<\frac{n+1}{2}$).

 Now, if $n$ is odd, then numerical computations indicate that $I(f_{Y_{[r:n]}},f_{Y})$ is increasing (decreasing) in $r$ for $1\leq r<\frac{n+1}{2}$,
 $0<\alpha\leq1$ ($\frac{n+1}{2}<r\leq n$, $-1\leq\alpha<0$).
\end{Example}
\begin{Example}
 Let $(X_i,Y_i)$, $i=1,2,\cdots,n$ be a random sample from Morgenstern type bivariate Logistic distribution with cdf
\begin{equation*}
F(x,y)=\left(1+\exp(-x)\right)^{-1}\left(1+\exp(-y)\right)^{-1}\left(1+\frac{\alpha
 e^{-x-y}}{(1+e^{-x})(1+e^{-y})}\right).
\end{equation*}
Computation shows that
\begin{equation}\label{eq20}
I(f_{Y_{[r:n]}},f_{Y})=1-0.6\alpha \left(\frac{n-2r+1}{n+1} \right).
\end{equation}

 By using \eqref{eq20}, we get
\begin{eqnarray*}
 D_{\alpha}(n)= I(f_{Y_{[n:n]}},f_{Y})- I(f_{Y_{[1:n]}},f_{Y})=1.2\alpha\left(\frac{n-1}{n+1}\right),
\end{eqnarray*}
which is positive, negative or zero whenever $0 <\alpha\leq1, n > 1$; $-1\leq \alpha < 0, n > 1$ or
$n = 1 \; or \; \alpha = 0$, respectively.

\end{Example}
\begin{Example}
Let $(X_{i},Y_{i})$, $i=1,2,...,n$ be a random sample from
Morgenstern type bivariate Rayleigh distribution with cdf
 \begin{equation*}
F(x,y)=\left(1-\exp(-\frac{x^{2}}{2\sigma_{1}^{2}})\right)\left(1-\exp(-\frac{y^{2}}{2\sigma_{2}^{2}})\right)\left(1+\alpha\exp\left(-\frac{x^{2}}{2\sigma_{1}^{2}}-\frac{y^{2}}{2\sigma_{2}^{2}}\right)\right).
\end{equation*}
From \eqref{Ig2}, we find
\begin{equation}
\label{eq3}
I(f_{Y_{[r:n]}},f_{Y})=\frac{\alpha(n-2r+1)}{n+1}(\log\sqrt{2}-\frac{1}{2})+1-\frac{1}{2}\psi(1)+\log(\frac{\sigma_{2}}{\sqrt{2}}).
\end{equation}
Using \eqref{eq3}, we have
\begin{eqnarray*}
 W_{\alpha}(n)= I(f_{Y_{[n:n]}},f_{Y})- I(f_{Y_{[1:n]}},f_{Y})=2\alpha\left(0.5-\log\sqrt{2}\right)\left(\frac{n-1}{n+1}\right),
\end{eqnarray*}
which is positive, negative or zero whenever $0 <\alpha\leq1, n > 1$; $-1\leq\alpha < 0, n > 1$ or
$n =1 \; or\; \alpha = 0$, respectively.
\end{Example}
\begin{Example}

Let $(X_{i},Y_{i})$, $i=1,2,...,n$ be a random sample from
Morgenstern type bivariate generalized exponential distribution (MTBGED) with cdf
 \begin{eqnarray*}
F_{X,Y}(x,y)=\{(1-e^{-\theta_{1}x})(1-e^{-\theta_{2}y})\}^{\lambda}[1+\alpha(1-(1-e^{-\theta_{1}x})^{\lambda})(1-(1-e^{-\theta_{2}y})^{\lambda})].
\end{eqnarray*}
By using \eqref{Ig2}, we get
 \begin{eqnarray}
 \label{eq4}
I(f_{Y_{[r:n]}},f_{Y})=-\log(\lambda\theta_2)
+B(\lambda)-\frac{\alpha(n-2r+1)}{n+1}D(\lambda)
+\frac{\lambda-1}{\lambda}[1+\frac{\frac{\alpha(n-2r+1)}{n+1}}{2}],
\end{eqnarray}
where $B(\lambda)=\psi \left(\lambda +1\right)-\psi \left(1\right)$ and $D(\lambda)=B(2\lambda)-B(\lambda)$. Using \eqref{eq4}, we have
\begin{eqnarray*}
 Q_{\alpha,\lambda}(n)= I(f_{Y_{[n:n]}},f_{Y})- I(f_{Y_{[1:n]}},f_{Y})=\frac{\alpha(n-1)}{n+1}\left[2D(\lambda)-\frac{\lambda-1}{\lambda}\right],
\end{eqnarray*}
which is positive, negative or zero whenever $0 <\alpha\leq1, n > 1$; $-1\leq\alpha < 0, n > 1$ or
$n = 1 \; or\; \alpha = 0$, respectively.
\end{Example}

\begin{Remark}
Let $(X_{i},Y_{i})$, $i=1,2,\cdots,n$ be a
random sample of size $n$ with pdf \eqref{Mf}. Then, from \eqref{Ig2} we have
\begin{equation*}
H(Y)=\frac{I(f_{Y_{[n:n]}},f_{Y})+I(f_{Y_{[1:n]}},f_{Y})}{2}.
\end{equation*}
\end{Remark}

\begin{Remark}
Let $(X_{i},Y_{i})$, $i=1,2,\cdots,n$ be a
random sample of size $n$ with pdf \eqref{Mf}. If $\lambda\geq1$ is an integer number and we change $r$ to $r\lambda$ and $n$ to $(n+1)\lambda-1$. Then, from \eqref{Ig2} we have
\begin{equation*}
I(f_{Y_{[r:n]}},f_{Y})=I(f_{Y_{[r\lambda:(n+1)\lambda-1]}},f_{Y}).
\end{equation*}
\end{Remark}

We consider the concomitants of order statistics whenever \\ $(X_1,Y_1),(X_2,Y_2),\ldots,(X_n,Y_n)$ are independent but otherwise arbitrarily distributed. Let us consider the Morgenstern family with cdf
\begin{equation}{\label{FGM}}
F_{X_i,Y_i}(x,y)=F_{X_i}(x)F_{Y_i}(y)\left[1+\alpha_i(1-F_{X_i}(x))(1-F_{Y_i}(y))\right]. \end{equation}
Now, suppose that $F_{X_i}(x)=F_{X}(x)$ , $F_{Y_i}(y)=F_{Y}(y)$ and $|\alpha_i|\leq 1$. Then in this particular case, the pdf’s of $Y_{[1:n]}$ and $Y_{[n:n]}$ are given by [6] as follows:
\begin{equation}
f_{[1:n]}(y)=f_Y(y)\left[1+\frac{n-1}{(n+1)n}\sum_{j=1}^n\alpha_j (1-2F_Y(y))\right],
\end{equation}
\begin{equation}
f_{[n:n]}(y)=f_Y(y)\left[1-\frac{n-1}{(n+1)n}\sum_{j=1}^n\alpha_j (1-2F_Y(y))\right].
\end{equation}

Now, in the following, the measures of inaccuracy for concomitants of extremes of order
statistics is represented.
\begin{Example}
Let $(X_i,Y_i),\, i=1,2,\ldots,n$ be independent random vectors from \eqref{FGM}. If $Y_{[1:n]}$ and $Y_{[n:n]}$ are concomitants of extremes of order statistics, then
\begin{equation}{\label{26}}
I(f_{[1:n]},f_{Y})=\left(1+\frac{n-1}{(n+1)n}\sum_{j=1}^n\alpha_j\right) H(Y) +2\frac{n-1}{(n+1)n}\sum_{j=1}^n\alpha_j\phi_f(u),
\end{equation}
\begin{equation}{\label{27}}
I(f_{[n:n]},f_{Y})=\left(1-\frac{n-1}{(n+1)n}\sum_{j=1}^n\alpha_j\right) H(Y) -2\frac{n-1}{(n+1)n}\sum_{j=1}^n\alpha_j\phi_f(u).
\end{equation}
By using \eqref{26} and \eqref{27} we have
\begin{equation*}
A_n=I(f_{[n:n]},f_{Y})-I(f_{[1:n]},f_{Y})=-\frac{2(n-1)}{n(n+1)}\Delta,
\end{equation*}
where we have set $\Delta=H(Y)\sum_{j=1}^n\alpha_j+2\sum_{j=1}^n\alpha_j\phi_f(u)$. If $\Delta>0$ ($\Delta<0$), then $A_n<0$ ($A_n>0$). Finally, we get

\begin{equation*}
I(f_{[n:n]},f_{Y})+I(f_{[1:n]},f_{Y})=2H(Y).
\end{equation*}

\end{Example}

\textbf{Case 2:} According to \eqref{Ig1}, if we put $m = -1$ and $k = 1$, then an inaccuracy measure between $f_{R_{[r]}}$ (density function of the concomitant of rth-record value) and $f_{Y}$ in Morgenstern family is obtained as follows:
\begin{eqnarray}\label{Ig3}
  I(f_{R_{[r]}},f_{Y})=\left(1+\alpha (2^{1-r}-1)\right) H(Y) +2\alpha (2^{1-r}-1)\phi_f(u).
\end{eqnarray}

\begin{Example}
Let $(X_{i},Y_{i})$, $i=1,2,...,n$ be a random sample of GBED with cdf
\begin{equation}
F(x,y)=\left(1-\exp\left(\frac{-x}{\theta_{1}}\right)\right)\left(1-\exp\left(\frac{-y}{\theta_{2}}\right)\right)\left[1+\alpha\exp\left(\frac{-x}{\theta_{1}}-\frac{y}{\theta_{2}}\right)\right].
\end{equation}
From \eqref{Ig3}, we find
\begin{eqnarray}
\label{eq22}
  I(f_{R_{[r]}},f_{Y})=[1+\log \theta_{2}]+\frac{\alpha}{2}\left(2^{1-r}-1\right).
\end{eqnarray}
By using \eqref{eq22}, we get

\begin{eqnarray*}
 A_{\alpha}(r)= I(f_{R_{[r]}},f_{Y})- I(f_{R_{[r-1]}},f_{Y})=-\alpha 2^{-r},
\end{eqnarray*}

which is positive, negative or zero whenever $-1\leq\alpha < 0, r > 1$; $0 <\alpha\leq1, r > 1$ or
$\alpha = 0$, respectively. Also, the difference between $ I(f_{R_{[r]}},f_{Y})$ and $H(Y)$ is
\begin{eqnarray*}
 B_{\alpha,n}(r)= I(f_{R_{[r]}},f_{Y})- H(Y)=\frac{\alpha}{2}\left(2^{1-r}-1\right).
\end{eqnarray*}
$B_{\alpha,n}(r)$ is positive, negative or zero whenever $-1\leq\alpha < 0, r > 1$; $0 <\alpha\leq1, r > 1$ or
$r = 1 \; or\; \alpha = 0$, respectively.

\end{Example}
\begin{Remark}
In analogy with \eqref{Ig1}, a measure of inaccuracy associated with $f_{Y}(y)$ and $g_{[r,n,m,k]}(y)$ is given  by
\begin{eqnarray*}
 I(f_{Y},g_{[r,n,m,k]})=H(Y)-E\left[\log \left(1+\alpha C^*(r,n,m,k)\left(1-2U\right)\right)\right],
\end{eqnarray*}
where U is uniformly distributed in (0,1).
\end{Remark}

Quantile functions are efficient alternatives to the distribution function in modelling and analysis of statistical data. The quantile function is defined by,
$$
 Q(u)=F^{-1}(u)=\inf\{y:F(y)\geq u\},\qquad  0< u< 1.
$$
Noting that $F(Q(u))=u$ and differentiating it with respect to $u$ yields \\
$q(u)f(Q(u))=1$. Let $Y$ be a nonnegative random variable with pdf $f(\cdot)$ and quantile function $Q(\cdot)$, then $f(Q(u))$ is called the density quantile function and $q(u) =Q'(u)$ is known as the quantile density function of $Y$. Now using \eqref{Ig1}, the corresponding quantile  based $I(g_{[r,n,m,k]},f_{Y})$ is defined as
\begin{eqnarray}
  I(g_{[r,n,m,k]},f_{Y})=E(\log q(U))+\alpha C^*(r,n,m,k)E\left[(1-2U)\log q(U)\right].
\end{eqnarray}

\subsection{CPI between $Y_{[r,n,m,k]}$ and $Y$}

If $Y_{[r,n,m,k]}$ is the concomitant of rth GOS from \eqref{Mf},
then the CPI measure between $G_{[r,n,m,k]}(y)$ and $F_{Y}(y)$ for $1\leq r \leq n$, $\alpha\neq 0$ is given by
\begin{eqnarray}{\label{IG1}}
  I(G_{Y_{[r,n,m,k]}},F_{Y})&=&-\int_0^{\infty}G_{[r,n,m,k]}(y)\log{F_{Y}(y)}dy \nonumber\\
  &=&\left[1+\alpha C^*(r,n,m,k)\right]{\mathcal{CE}}(Y)\nonumber\\
  &+&\alpha C^*(r,n,m,k)\int_0^{\infty}F^{2}_{Y}(y)\log{F_{Y}(y)}dy\nonumber\\
  &=&\left[1+\alpha C^*(r,n,m,k)\right]{\mathcal{CE}}(Y)\nonumber\\
  &-&\frac{\alpha}{2} C^*(r,n,m,k){\mathcal{CE}}(Y_{(2:2)}),
\end{eqnarray}
where ${\mathcal{CE}}(Y)$ and ${\mathcal{CE}}(Y_{(2:2)})$ are the cumulative entropy of the random variables $Y$ and $Y_{(2:2)}$, respectively (see [4]).

\begin{Remark}
In analogy with \eqref{IG1}, a measure of inaccuracy associated with $F_{Y}$ and $G_{[r,n,m,k]}$ is given by
\begin{eqnarray*}
 I(F_{Y},G_{[r,n,m,k]})=\mathcal{CE}(Y)-E\left[\frac{U\log{\left(1+\alpha C^*(r,n,m,k)\left(1-U\right)\right)}}{f(F^{-1}(U))}\right].
\end{eqnarray*}
\end{Remark}

\textbf{Case 1:} If we put $m = 0$ and $k = 1$, then a measure of inaccuracy between $F_{Y_{[r:n]}}$ (distribution function of rth concomitant of order statistic) and $F_{Y}$ is presented as
\begin{eqnarray}
  I(F_{Y_{[r:n]}},F_{Y})&=&-\int_0^{\infty}F_{Y_{[r:n]}}(y)\log{F_{Y}(y)}dy \nonumber\\
  &=&\left[1+\alpha\left(\frac{n-2r+1}{n+1}\right)\right]{\mathcal{CE}}(Y)\nonumber\\
  &+&\alpha\left(\frac{n-2r+1}{n+1}\right)\int_0^{\infty}F^{2}_{Y}(y)\log{F_{Y}(y)}dy\nonumber\\
  &=&\left[1+\alpha\left(\frac{n-2r+1}{n+1}\right)\right]{\mathcal{CE}}(Y)\nonumber\\
  &-&\alpha\left(\frac{n-2r+1}{2(n+1)}\right){\mathcal{CE}}(Y_{(2:2)}).
\end{eqnarray}

\begin{Example}
 Let $(X_i,Y_i)$, $i=1,2,\cdots,n$ be a random sample from Morgenstern type bivariate uniform distribution (MTBUD) with cdf
\begin{equation*}
F(x,y)=\frac{xy}{\theta_{1}\theta_{2}}\left[1+\alpha(1-\frac{x}{\theta_{1}})(1-\frac{y}{\theta_{2}})\right],\;\;0<x<\theta_{1},\;\;0<y<\theta_{2}.
\end{equation*}
Computation shows that
 \begin{eqnarray}
 \label{eq5}
  I(F_{Y_{[r:n]}},F_{Y})
  &=&\left[1+\alpha\left(\frac{n-2r+1}{n+1}\right)\right]\frac{\theta_{2}}{4}-\alpha\left(\frac{n-2r+1}{n+1}\right)\frac{\theta_{2}}{9}
  \nonumber\\
  &=&\frac{\theta_{2}}{4}+\alpha\left(\frac{n-2r+1}{n+1}\right)\frac{5\theta_{2}}{36}.
\end{eqnarray}
Using \eqref{eq5}, we have
\begin{eqnarray*}
 D_{\alpha,\theta_{2}}(n)= I(F_{Y_{[n:n]}},F_{Y})- I(F_{Y_{[1:n]}},F_{Y})=\frac{5\alpha\theta_{2}(-n+1)}{18(n+1)}.
\end{eqnarray*}
which is positive, negative or zero whenever $-1 \leq\alpha<0$, $0<\alpha \leq 1$ or $\alpha = 0$, respectively.

\end{Example}
\begin{Example}
 Let $(X_i,Y_i)$, $i=1,2,\cdots,n$ be a random sample from GBED. Then, computation shows that
 \begin{eqnarray}
 \label{eq6}
  I(F_{Y_{[r:n]}},F_{Y})&=&\left[1+\alpha\left(\frac{n-2r+1}{n+1}\right)\right]\left[\frac{\pi^{2}}{6}-1\right]\theta_{2}\nonumber\\
  &-&\alpha\left(\frac{n-2r+1}{n+1}\right)\left[\frac{\pi^{2}}{6}-\frac{5}{4}\right]\theta_{2}\nonumber\\
  &=&\left[\frac{\pi^{2}}{6}-1\right]\theta_{2}+\frac{\alpha\theta_{2}}{4}\left(\frac{n-2r+1}{n+1}\right).
  \end{eqnarray}
Using \eqref{eq6}, we have
\begin{eqnarray*}
 Q_{\alpha,\theta_{2}}(n)= I(F_{Y_{[n:n]}},F_{Y})- I(F_{Y_{[1:n]}},F_{Y})=\frac{\alpha\theta_{2}(-n+1)}{2(n+1)},
\end{eqnarray*}
which is positive, negative or zero whenever $-1\leq\alpha<0$, $0<\alpha \leq 1$ or $\alpha = 0$, respectively.
\end{Example}
\begin{Example}
 Let $(X_i,Y_i)$, $i=1,2,\cdots,n$ be a random sample from Morgenstern type bivariate inverse Weibull distribution  with cdf
\begin{eqnarray*}
F(x,y)&=&\exp\left[-\left(\frac{\theta_{1}}{x}\right)^{\beta_{1}}-\left(\frac{\theta_{2}}{y}\right)^{\beta_{2}}\right]\nonumber\\ &\times&\left[1+\alpha\left(1-\exp\left[-\left(\frac{\theta_{1}}{x}\right)^{\beta_{1}}\right]\right)\left(1-\exp\left[-\left(\frac{\theta_{2}}{y}\right)^{\beta_{2}}\right]\right)\right].
\end{eqnarray*}
Computation shows that
 \begin{eqnarray}
 \label{eq7}
  I(F_{Y_{[r:n]}},F_{Y})
  &=&\left[1+\alpha\left(\frac{n-2r+1}{n+1}\right)\right]\frac{\theta_{2}}{\beta_{2}}\Gamma\left(\frac{\beta_{2}-1}{\beta_{2}}\right)\nonumber\\
  &-&\alpha\left(\frac{n-2r+1}{n+1}\right)\frac{2^{\frac{1}{\beta_{2}}-1}\theta_{2}}{\beta_{2}}\Gamma\left(\frac{\beta_{2}-1}{\beta_{2}}\right)
  \nonumber\\
  &=&\frac{\theta_{2}}{\beta_{2}}\Gamma\left(\frac{\beta_{2}-1}{\beta_{2}}\right)+\alpha\left(\frac{n-2r+1}{n+1}\right)\frac{\theta_{2}}{\beta_{2}}\Gamma\left(\frac{\beta_{2}-1}{\beta_{2}}\right)\left(1-2^{\frac{1}{\beta_{2}}-1}\right).\nonumber\\
  \end{eqnarray}
Using \eqref{eq7}, we have
\begin{eqnarray*}
 D_{\alpha,\theta_{2}}(n)= I(F_{Y_{[n:n]}},F_{Y})- I(F_{Y_{[1:n]}},F_{Y})=\alpha\frac{\theta_{2}(1-n)}{\beta_{2}(n+1)}\Gamma\left(\frac{\beta_{2}-1}{\beta_{2}}\right)\left(1-2^{\frac{1}{\beta_{2}}-1}\right).
\end{eqnarray*}
which is positive, negative or zero whenever $-1 \leq\alpha<0$, $0<\alpha \leq 1$ or $\alpha = 0$, respectively.
\end{Example}

\begin{Proposition}
Let $(X_i,Y_i)$, $i=1,2,\cdots,n$ be a random sample from Morgenstern family. Then for $1\leq r\leq\frac{n+1}{2}$, we have
\begin{equation}
  I(F_{Y_{[r:n]}},F_{Y})\leq(\geq){\mathcal{CE}}(Y), \; \; -1\leq\alpha<0 \;(0<\alpha\leq1).
\end{equation}

{\bf Proof.} The proof follows by recalling Proposition 4.8 of [4] .
\end{Proposition}

\textbf{Case 2:} If we put $m = -1$ and $k = 1$, then a measure of inaccuracy between $F_{R_{[r]}}$ (distribution function of nth concomitant of upper record value) and $F_{Y}$ is presented as

\begin{eqnarray}
 I(F_{R_{[r]}},F_{Y})&=&-\int_0^{\infty}F_{R_{[r]}}(y)\log{F_{Y}(y)}dy \nonumber\\
  &=&[1+\alpha(2^{1-r}-1)]{\mathcal{CE}}(Y)+\alpha(2^{1-r}-1)\int_{0}^{\infty}F^{2}_{Y}(y)\log F_{Y}(y)dy \nonumber\\
  &=&[1+\alpha(2^{1-r}-1)]{\mathcal{CE}}(Y)-\frac{\alpha}{2}(2^{1-r}-1){\mathcal{CE}}(Y_{(2:2)}).
    \end{eqnarray}
\begin{Proposition}
Let $(X_i,Y_i)$, $i=1,2,\cdots,n$ be a random sample from Morgenstern family. Then, we have
\begin{equation}
  I(F_{R_{[r]}},F_{Y})\leq(\geq){\mathcal{CE}}(Y), \; \; 0<\alpha\leq1(-1\leq\alpha<0).
\end{equation}

{\bf Proof.} The proof follows by recalling  Proposition 4.8 of [4] .
\end{Proposition}

\section{Empirical CPI for concomitants of GOS}
In this section we address the problem of estimating the CPI for concomitants by means of the empirical CPI. Let $(X_{i},Y_{i})$, $i=1,2,\cdots,n$ be a
random sample of size $n$ from Morgenstern family. Then according to \eqref{IG1}, the empirical CPI between $G_{Y_{[r,n,m,k]}}$  and $F_{Y}$ can be obtained as follows:
\begin{eqnarray}{\label{Ihat1}}
  \widehat{I}(G_{Y_{[r,n,m,k]}},F_{Y})&=&\left[1+\alpha C^*(r,n,m,k)\right]\sum_{j=1}^{n-1}U_{j}\left(\frac{j}{n}\right)\left(-\log \frac{j}{n}\right)\nonumber\\
  &-&\alpha C^*(r,n,m,k)\sum_{j=1}^{n-1}U_{j}\left(\frac{j}{n}\right)^2\left(-\log \frac{j}{n}\right)\nonumber\\
  &=&\sum_{j=1}^{n-1}U_{j}\left(\frac{j}{n}\right)\left(-\log \frac{j}{n}\right)\left[1+\alpha C^*(r,n,m,k)\left(1-\frac{j}{n}\right) \right],\nonumber\\
\end{eqnarray}
where $U_{j}=Z_{(j+1)}-Z_{(j)}, j=1,2,...,n-1$ are the sample spacings based on ordered random samples.

\textbf{Case 1:} If we put $m = 0$ and $k = 1$, then the empirical CPI between $F_{Y_{[r:n]}}$ and $F_{Y}$ is given by
\begin{eqnarray}
\widehat{I}(F_{Y_{[r:n]}},F_{Y})
=\sum_{j=1}^{n-1}U_{j} \frac{j}{n}\left(-\log \frac{j}{n}\right)\left[1+\alpha\left( \frac{n-2r+1}{n+1}\right)\left(1-\frac{j}{n}\right) \right].
\end{eqnarray}

\textbf{Case 2:} If we put $m = -1$ and $k = 1$, then the empirical CPI between $F_{R_{[r]}}$ and $F_{Y}$ can be written as
\begin{eqnarray}{\label{Ihat4}}
\widehat{I}(F_{R_{[r]}},F_{Y})
=\sum_{j=1}^{n-1}U_{j} \frac{j}{n}\left(-\log \frac{j}{n}\right)\left[1+\alpha\left(2^{1-r}-1\right)\left(1-\frac{j}{n}\right) \right].
\end{eqnarray}

\begin{Example}
Let $(X_{i},Y_{i})$, $i=1,2,...,n$ be a random sample from MTBGED with $\lambda=1$, then the sample spacings $U_{j}$ are independent and exponentially
distributed with mean $\frac{1}{\theta_2 (n-j)}$ (for more details see [10]). Now from \eqref{Ihat4} we obtain
\begin{eqnarray}
E[\widehat{I}(F_{R_{[r]}},F_{Y})]=\frac{1}{\theta_2 }\sum_{j=1}^{n-1}\frac{j}{n(n-j)}\left(-\log \frac{j}{n}\right)\left[1+\alpha\left(2^{1-r}-1\right)\left(1-\frac{j}{n}\right)\right],
\end{eqnarray}
and
\begin{eqnarray}
Var[\widehat{I}(F_{R_{[r]}},F_{Y})]=\frac{1}{\theta_2 ^2}\sum_{j=1}^{n-1}\left(\frac{j}{n(n-j)}(-\log \frac{j}{n})\left[1+\alpha\left(2^{1-r}-1\right)\left(1-\frac{j}{n}\right)\right]\right)^2.\nonumber\\
\end{eqnarray}
We have computed the values of $E[\widehat{I}(F_{R_{[r]}},F_{Y})]$ and $Var[\widehat{I}(F_{R_{[r]}},F_{Y})]$ for sample sizes $n=10, 15,20$, $\theta_2 =0.5,1,2$, $\alpha=-1,-0.5,0.5,1$ and $r=2$ in Table 1. We can easily see that $E[\widehat{I}(F_{R_{[r]}},F_{Y})]$ and $Var[\widehat{I}(F_{R_{[r]}},F_{Y})]$ are decreasing in $\alpha$  and  $\theta_{2}$. Also, we consider that $\lim_{n\rightarrow\infty}Var[\widehat{I}(F_{R_{[r]}},F_{Y})]=0$.
\end{Example}

\begin{Example}
 Let $(X_i,Y_i)$, $i=1,2,\cdots,n$ be a random sample from \\
  MTBUD with $\theta_1=\theta_2=1$. Then the sample spacings $U_{j}$ are independent of beta distribution with parameters 1 and $n$ (for more details see [10]). Now from \eqref{Ihat4} we
obtain
\begin{eqnarray}
E[\widehat{I}(F_{R_{[r]}},F_{Y})]=
\frac{1}{n+1}\sum_{j=1}^{n-1}\frac{j}{n}\left(-\log \frac{j}{n}\right)\left[1+\alpha\left(2^{1-r}-1\right)\left(1-\frac{j}{n}\right)\right],
\end{eqnarray}
and
\begin{eqnarray}
Var[\widehat{I}(F_{R_{[r]}},F_{Y})]=
\frac{n}{(n+1)^2(n+2)}\sum_{j=1}^{n-1}\left(\frac{j}{n}(-\log \frac{j}{n})\left[1+\alpha\left(2^{1-r}-1\right)\left(1-\frac{j}{n}\right)\right]\right)^2.\nonumber\\
\end{eqnarray}
We have computed the values of $E[\widehat{I}(F_{R_{[r]}},F_{Y})]$ and $Var[\widehat{I}(F_{R_{[r]}},F_{Y})]$ for sample sizes $n=10, 15,20$, $\alpha=-1,-0.5,0.5,1$ and $r=2$ in Table 2. We can easily see that $E[\hat{I}(F_{R_{[r]}},F_{Y})]$  and  $Var[\widehat{I}(F_{R_{[r]}},F_{Y})]$  are decreasing in $\alpha$. Also, we consider that $\lim_{n\rightarrow\infty}Var[\widehat{I}(F_{R_{[r]}},F_{Y})]=0$.
\end{Example}

\begin{table}[h]
\begin{center}
\caption{Numerical values  of $\mathbb{E}[\widehat{I}(F_{R_{[r]}},F_{Y})]$ and $Var[\widehat{I}(F_{R_{[r]}},F_{Y})]$ for MTBGED with $\lambda=1$.}\label{tab.1}
{\small
\begin{tabular}{|c|ccc|ccc|ccc|ccc|} \hline
 &   \multicolumn{12}{c|}{$\mathbb{E}[\widehat{I}(F_{R_{[r]}},F_{Y})] $}   \\ \hline
$\theta_2 $  & 0.5 & 1 & 2 & 0.5 & 1 & 2 & 0.5 & 1 & 2& 0.5 & 1 & 2 \\   \hline
$n$ &  & $\alpha=-1$ &  &        & $\alpha=-0.5$ &  &            & $\alpha=0.5$  &  &          & $\alpha=1$ & \\   \hline
10&   1.429 &  0.714   &  0.357  &  1.306  &  0.653  &  0.326  & 1.061 &  0.530  &  0.265  & 0.938 &  0.469  &  0.234  \\
15&   1.468  &  0.734   &  0.367  & 1.344  &  0.672  &  0.336  & 1.096 &  0.548  &  0.274  & 0.972 &  0.486  &  0.243  \\
20&   1.487  &  0.743   &  0.372  & 1.362  &  0.681  &  0.340  & 1.114 &  0.557  &  0.278  & 0.989 &  0.494  &  0.247  \\

  \hline \hline
  &   \multicolumn{12}{c|}{$Var[\widehat{I}(F_{R_{[r]}},F_{Y})]$}   \\ \hline
$\theta_2 $  & 0.5 & 1 & 2 & 0.5 & 1 & 2 & 0.5 & 1 & 2& 0.5 & 1 & 2 \\   \hline
$n$ &  & $\alpha=-1$ &  &    & $\alpha=-0.5$ &  &      & $\alpha=0.5$  &  &     & $\alpha=1$ & \\   \hline
10  & 0.241  & 0.060   & 0.015  & 0.205  & 0.051   &  0.013 & 0.144  & 0.036  &  0.009  & 0.119  & 0.030  & 0.007  \\
15  & 0.165  &  0.041  &  0.010 & 0.141  &  0.035  & 0.009  & 0.100  & 0.025  & 0.006  & 0.083  & 0.021  & 0.005  \\
20  & 0.126  &  0.031  &  0.008 & 0.108  & 0.027   & 0.007  & 0.077  & 0.019  & 0.005  &  0.064 &  0.016 &  0.004 \\
  \hline
\end{tabular}
}
\end{center}
\end{table}

\begin{table}[h]
\begin{center}
\caption{Numerical values of $\mathbb{E}[\widehat{I}(F_{R_{[r]}},F_{Y})]$ and $Var[\widehat{I}(F_{R_{[r]}},F_{Y})]$ for MTBUD with $\theta_1=\theta_2=1$.} \label{tab.1}
\begin{tabular}{|c|cccc|cccc|}
\hline
    &  \multicolumn{4}{|c}{ $\mathbb{E}[\widehat{I}(F_{R_{[r]}},F_{Y})]$ }               & \multicolumn{4}{|c|}{  $Var[\widehat{I}(F_{R_{[r]}},F_{Y})]$  }          \\ \hline

    $n$   & $\alpha=-1$ & $\alpha=-0.5$  & $\alpha=0.5$ & $\alpha=1$     & $\alpha=-1$ & $\alpha=-0.5$ & $\alpha=0.5$ & $\alpha=1$\\ \hline
      10   & 0.285 & 0.254 & 0.192 & 0.162    & 0.008  &  0.007 & 0.004  & 0.003  \\
      15   & 0.297 & 0.264 &0.200  &  0.168   & 0.006  &  0.005 &  0.003 &  0.002 \\
      20   & 0.302 & 0.270 & 0.204 & 0.171    &  0.005 &  0.004 & 0.002  & 0.001  \\
                                         \hline
\end{tabular}
\end{center}
\end{table}

 \begin{Theorem}
Let $(X_{i},Y_{i})$, $i=1,2,\cdots,n$ be a
random sample of size $n$ from Morgenstern family. Then we have
\begin{eqnarray*}
\widehat{I}(F_{R_{[r]}},F_{Y})\,\longrightarrow I(F_{R_{[r]}},F_{Y})\; \; \; a.s\; as \; n\rightarrow \infty.
\end{eqnarray*}

{\bf Proof.} From relation \eqref{Ihat4}, we obtain
\begin{eqnarray*}
 \widehat{I}(F_{R_{[r]}},F_{Y})=\left[1+\alpha\left(2^{1-r}-1\right)\right]\widehat{\mathcal{CE}}(Y) -\frac{\alpha}{2}\left(2^{1-r}-1\right)\widehat{\mathcal{CE}}(Y_{(2:2)}).
\end{eqnarray*}

 Since $\widehat{\mathcal{CE}}(Y)\longrightarrow \mathcal{CE}(Y)$ and $\widehat{\mathcal{CE}}(Y_{(2:2)})\longrightarrow \mathcal{CE}(Y_{(2:2)})$, then proof follows by [4].

\end{Theorem}

\begin{Theorem}
Let $(X_{i},Y_{i})$, $i=1,2,...,n$ be a random sample from MTBGED with $\lambda=1$, then
\begin{eqnarray*}
Z_n:= \frac{\widehat{I}(F_{R_{[r]}},F_{Y})-E\left[\widehat{I}(F_{R_{[r]}},F_{Y})\right]}{\sqrt{Var\left[\widehat{I}(F_{R_{[r]}},F_{Y})\right]}}
\end{eqnarray*}
converges in distribution to a standard normal variable as $n\rightarrow \infty$.

{\bf Proof.} First the empirical measure $\hat{I}(F_{R_{[r]}},F_{Y})$ can be expressed as the following sum of
independent random variables as
\begin{eqnarray*}
\widehat{I}(F_{R_{[r]}},F_{Y})=\sum_{j=1}^{n-1} W_j,
\end{eqnarray*}
where  $W_j=U_{j}\frac{j}{n}\left(-\log \frac{j}{n}\right)\left[1+\alpha\left(2^{1-r}-1\right)\left(1-\frac{j}{n}\right)\right]$ are independent random variables with the
mean and variance given by
\begin{equation*}
E[W_j]=\frac{1}{n\theta_2 (1-\frac{1}{j/n})}\left(\log\frac{j}{n}\right)\left[1+\alpha\left(2^{1-r}-1\right)\left(1-\frac{j}{n}\right)\right],
\end{equation*}
\begin{equation*}
Var[W_j]=\frac{1}{n^{2}\theta^{2}_2 (1-\frac{1}{j/n})^{2}}\left(\log\frac{j}{n}\right)^{2}\left[1+\alpha\left(2^{1-r}-1\right)\left(1-\frac{j}{n}\right)\right]^{2}.
\end{equation*}

Since $E[|W_j-E(W_j)|^3]=2e^{-1}(6-e)[E(W_j)]^3$ for any exponentially distributed random
variable $W_j$, by setting $\alpha_{j,k}=E[|W_j-E(W_j)|^k]$ the following approximations hold for large n:
\begin{eqnarray*}
\sum_{j=1}^n \alpha_{j,2}&=&\frac{1}{n^2\theta_2^2}\sum_{j=1}^n\frac{1}{\left(1-\frac{1}{j/n}\right)^2}\left(\log\frac{j}{n}\right)^2\left[1+\alpha\left(2^{1-r}-1\right)\left(1-\frac{j}{n}\right)\right]^2 \nonumber \\
 &\approx& \frac{c_2}{n\theta_2^2},\\
\sum_{j=1}^n \alpha_{j,3}&=&\frac{2(6-e)}{en^3\theta_2^3}\sum_{j=1}^n\frac{1}{\left(1-\frac{1}{j/n}\right)^3}\left(\log\frac{j}{n}\right)^3\left[1+\alpha\left(2^{1-r}-1\right)\left(1-\frac{j}{n}\right)\right]^3\nonumber \\
 &\approx& \frac{2(6-e)c_3}{en^2\theta_2^3},
\end{eqnarray*}
where
\begin{eqnarray*}
c_k:=\int_0^1 \left(\frac{\log x}{1-1/x}\right)^k\left[1+\alpha\left(2^{1-r}-1\right)\left(1-x\right)\right]^k.
\end{eqnarray*}
Hence, Lyapunov’s condition of the central limit theorem is satisfied (see
[3]):
\begin{eqnarray*}
\frac{(\alpha_{1,3}+\cdots+\alpha_{n,3})^{1/3}}{(\alpha_{1,2}+\cdots+\alpha_{n,2})^{1/2}}\approx \frac{[2(6-e)c_3]^{1/3}}{e^{1/3}c_2^{1/2}}n^{-1/6} \rightarrow 0 \;\;\;\;\; as \;\; n\rightarrow \infty,
\end{eqnarray*}
which completes the proof.

\end{Theorem}

\section*{Acknowledgement}
The authors are thankful to the referee and editor for their
valuable suggestions towards the improvement of the paper.

\end{article}
\end{document}